\documentclass[twocolumn,aps,nofootinbib,secnumarabic,floatfix]{revtex4}
\usepackage{amssymb,amsmath2000,amsthm2000,times,mathptmx,pstricks,pst-node}
\newtheorem{theorem}{Theorem}[section]
\newtheorem{proposition}[theorem]{Proposition}
\newtheorem{lemma}[theorem]{Lemma}
\theoremstyle{remark}
\newtheorem*{remark}{Remark}
\newcommand{\bqn}{\begin{equation}}
\newcommand{\eqn}{\end{equation}}
\newcommand{\C}{\mathbb{C}}
\newcommand{\Z}{\mathbb{Z}}
\newcommand{\Q}{\mathbb{Q}}
\newcommand{\ahat}{\hat{a}}
\newcommand{\bhat}{\hat{b}}
\newcommand{\chat}{\hat{c}}
\newcommand{\Det}{\mathrm{Det}}
\newcommand{\End}{\mathrm{End}}
\newcommand{\ch}{\mathrm{ch}}
\renewcommand{\tensor}{\otimes}
\renewcommand{\sl}{\mathrm{sl}}
\newcommand{\ie}{\emph{i.e.}}

\newenvironment{fullfigure}[2]
    {\begin{figure}[ht]\begin{center}\def\ffa{#1}\def\ffb{#2}}
    {\vspace{\baselineskip}\caption{\ffb.}\label{\ffa}\end{center}\end{figure}}

\psset{linewidth=.25pt,dash=3pt 3pt}
\newrgbcolor{red1}{1 .25 .25}
\newrgbcolor{red2}{1 .5 .5}
\newrgbcolor{red3}{1 .75 .75}
\newcommand{\whitedot}{\pscircle[fillstyle=solid,fillcolor=white](0,0){2pt}}
\newcommand{\blackdot}{\qdisk(0,0){2pt}}
\newcommand{\td}[3]{%
! #1 -.6 mul #2 .8 mul add #1 -4 mul #2 -3 mul add #3 12 mul add 13 div}

\begin{document}

\title{Four symmetry classes of plane partitions under one roof}
\author{Greg Kuperberg}
\affiliation{Department of Mathematics, Yale University, New Haven, CT 06520}
\thanks{Supported by an NSF Postdoctoral Fellowship, grant \#DMS-9107908}
\email[Current email: ]{greg@math.ucdavis.edu}


\begin{abstract}
In previous paper, the author applied the permanent-determinant method of
Kasteleyn and its non-bipartite generalization, the Hafnian-Pfaffian method, to
obtain a determinant or  a Pfaffian that enumerates each of the ten symmetry
classes of plane partitions.  After a cosmetic generalization of the Kasteleyn
method, we identify the matrices in the four determinantal cases (plain plane
partitions, cyclically symmetric plane partitions, transpose-complement plane
partitions, and the intersection of the last two types) in the representation
theory of $\sl(2,\C)$. The result is a unified proof of the four enumerations.
\end{abstract}

\maketitle

Stanley \cite{Stanley:symmetries} and Robbins recognized that there are
ten symmetry classes of plane partitions in a box, and, with the aid of
computer experiments and other work, that the number of plane
partitions in each symmetry class is given by a product formula. The
program of proving each of these formulas has recently been completed
\cite{Andrews:tsscpp,Kuperberg:perdet,Stembridge:enumeration}, but
although the formulas are very similar, there is no known unified
treatment.  Moreover, some of the enumerations presently require
difficult, ad-hoc calculations with generating functions or matrices.
In fact, the enumeration of one of the classes (totally symmetric plane
partitions or TSPP's) admits a natural $q$-analogue which is still
open.

In a previous paper, the author \cite{Kuperberg:perdet} discussed the
permanent-determinant method of Kasteleyn and its generalization, the
Hafnian-Pfaffian method, as the first step towards a possible unified
treatment.  (The connection between this method and plane partitions was
discovered jointly with James Propp.)  In this paper, we carry out this
strategy for the four symmetry class that are given by determinants and not
merely Pfaffians: Unrestricted plane partitions (PP's), cyclically symmetric
plane partitions (CSPP's), transpose-complement plane partitions (TCPP's),
and cyclically symmetric, transpose-complement plane partitions (CSTCPP's). 
The new idea is to find the matrices given by the permanent-determinant
method in the representation theory of the Lie algebra $\sl(2,\C)$. The main
result of the present paper is the following theorem:

\begin{theorem} A Kasteleyn-flat, weighted adjacency matrix for the graph
$Z(a,b,c)$, whose matchings are bijective with plane partitions in an $a
\times b \times c$ box, arises as $\alpha(X)|_{-1}$ in a representation
$\alpha$ of $\sl(2,\C)$ which is a tensor product of three irreducible
representations, with a basis formed from weight bases of irreducible
representations.  Consequently, $(\det \alpha(X)|_{-1})/m$
is the number of such plane partitions, where $m$ is the value of 
any term in the determinant.
\end{theorem}

Interestingly, the representation theory of the quantum group $U_q(\sl(2,\C))$
similarly yields the $q$-enumeration of PP's, but the author could not obtain
the known $q$-enumeration of CSPP's by this method.

Proctor's minuscule method is a competing method that, in various
forms, yields the enumeration of five symmetry classes of plane
partitions
\cite{Kuperberg:minuscule,Proctor:bruhat,Stembridge:minuscule}. Proctor's
method partly explains Stembridge's $q=-1$ phenomenon
\cite{Stembridge:q-1}, and it is applicable to any symmetry class not
involving cyclic symmetry. Likewise, the method given here is
applicable to any symmetry class which does not involve transposition or
complementation, but which may involve their product.

The author believes that one of these two methods, or perhaps some
combination of the two, should lead to a complete unified enumeration.
Both methods involve Cartan-Weyl representation theory.  Moreover, many
of the enumerations can be analyzed with the Gessel-Viennot
non-intersecting lattice-path method \cite{GV:binomial}, which in turn
can be recognized as a condensed version of the permanent-determinant
method \cite{Kuperberg:perdet}. For example, the Gessel-Viennot method
applied to PP's yield Carlitz matrices.  The determinants of these
matrices are well-known; interestingly, Proctor
\cite{Proctor:lefschetz} rederived them, and therefore enumerated PP's,
by an argument which is parallel to the one presented here.

The author would like to thank James Propp for suggesting numerous
corrections and for his continued interest in this work.

\section{Preliminaries}

\subsection{Plane partitions and their symmetries}

A \emph{ plane partition} in an $a \times b \times c$ box is a collection
of unit cubes in the rectangular solid $[0,a] \times [0,b] \times
[0,c]$ which is stable under gravitational attraction towards the
origin. Equivalently, it is an $a \times b$ matrix of integers with
entries between 0 and $c$; the stack of cubes is then a \emph{ Ferrer's
diagram} (or bar graph) of the matrix.  Figure~\ref{ftiling} shows a
picture of a plane partition together with three back walls of the
box.  We let $N(a,b,c)$ be the total number of plane partitions in an
$a \times b \times c$ box.

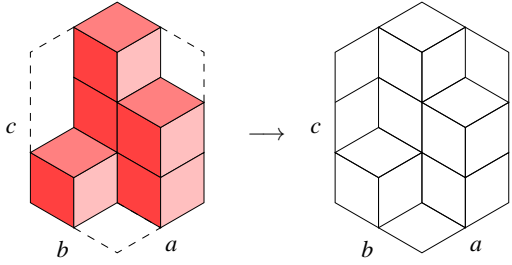
\begin{fullfigure}{ftiling}{A plane partition in a box and its lozenge tiling}
\psset{xunit=.577cm,yunit=.667cm}
\pspicture[.45](-3,0)(3,5)
\psset{fillstyle=solid}
\pspolygon[fillcolor=red1](-1,4.5)(-1,3.5)(0,3)(0,4)
\pspolygon[fillcolor=red3](1,4.5)(1,3.5)(0,3)(0,4)
\pspolygon[fillcolor=red1](-1,3.5)(-1,2.5)(0,2)(0,3)
\pspolygon[fillcolor=red1](0,3)(0,2)(1,1.5)(1,2.5)
\pspolygon[fillcolor=red3](2,3)(2,2)(1,1.5)(1,2.5)
\pspolygon[fillcolor=red1](0,2)(0,1)(1,.5)(1,1.5)
\pspolygon[fillcolor=red3](2,2)(2,1)(1,.5)(1,1.5)
\pspolygon[fillcolor=red1](-2,2)(-2,1)(-1,.5)(-1,1.5)
\pspolygon[fillcolor=red3](0,2)(0,1)(-1,.5)(-1,1.5)
\pspolygon[fillcolor=red2](-1,4.5)(0,4)(1,4.5)(0,5)
\pspolygon[fillcolor=red2](0,3)(1,2.5)(2,3)(1,3.5)
\pspolygon[fillcolor=red2](-2,2)(-1,1.5)(0,2)(-1,2.5)
\psset{fillstyle=none}
\psline[linestyle=dashed](-2,2)(-2,4)(-1,4.5)
\psline[linestyle=dashed](1,4.5)(2,4)(2,3)
\psline[linestyle=dashed](-1,.5)(0,0)(1,.5)(0,1)
\pcline[linestyle=none](-2,1)(-2,4)\Aput{$c$}
\pcline[linestyle=none](0,0)(-2,1)\Aput{$b$}
\pcline[linestyle=none](0,0)(2,1)\Bput{$a$}
\endpspicture $\longrightarrow$
\pspicture[.45](-3,0)(3,5)
\pspolygon(-1,4.5)(-1,3.5)(0,3)(0,4)
\pspolygon(1,4.5)(1,3.5)(0,3)(0,4)
\pspolygon(-1,3.5)(-1,2.5)(0,2)(0,3)
\pspolygon(0,3)(0,2)(1,1.5)(1,2.5)
\pspolygon(2,3)(2,2)(1,1.5)(1,2.5)
\pspolygon(0,2)(0,1)(1,.5)(1,1.5)
\pspolygon(2,2)(2,1)(1,.5)(1,1.5)
\pspolygon(-2,2)(-2,1)(-1,.5)(-1,1.5)
\pspolygon(0,2)(0,1)(-1,.5)(-1,1.5)
\pspolygon(-1,4.5)(0,4)(1,4.5)(0,5)
\pspolygon(0,3)(1,2.5)(2,3)(1,3.5)
\pspolygon(-2,2)(-1,1.5)(0,2)(-1,2.5)
\psline(-2,2)(-2,4)(-1,4.5)
\psline(1,4.5)(2,4)(2,3)
\psline(-1,.5)(0,0)(1,.5)(0,1)
\psline(-2,3)(-1,3.5)
\pcline[linestyle=none](-2,1)(-2,4)\Aput{$c$}
\pcline[linestyle=none](0,0)(-2,1)\Aput{$b$}
\pcline[linestyle=none](0,0)(2,1)\Bput{$a$}
\endpspicture
\end{fullfigure}

There are three natural symmetry operations on plane partitions:

\begin{description}
\item[1.]  \emph{ Transposition}, $\tau$, is a reflection through a
diagonal plane of the box of a plane partition.  A
transposition-symmetric plane partition is called, briefly, a
\emph{ symmetric} plane partition.
\item[2.]  \emph{ Rotation}, $\rho$, is a rotation by $120^\circ$ about
the long diagonal of the box of a plane partition.  Symmetry under
rotation is also called \emph{ cyclic symmetry}.
\item[3.]  \emph{ Complementation}, $\kappa$, takes the point
$(x,y,z)$ to the point $(a-x,b-y,c-z)$ for an $a \times b
\times c$ box.  The operation $\kappa$ is the same for individual
cubes, but for a plane partition $P$, $\kappa(P)$ is the set of
all cubes $C$ such that $\kappa(C) \notin P$.
\end{description}

By convention, the operations $\tau$ and $\rho$ are only defined when
the hexagon or box possesses them as symmetries, while the operation
$\kappa$ always exists.  The group of all symmetries has ten conjugacy
classes of subgroups, and for a subgroup $G$, let $N_G(a,b,c)$ be the
number of plane partitions with that symmetry.  We describe $G$ by
generators as $G = \langle g_1, g_2,\ldots \rangle$ and we often omit
the angle brackets, so that $N_{\rho,\tau}(a,a,a)$ is the number of
totally symmetric plane partitions.  (By convention, the phrase ``total
symmetry'' does not include the symmetry operation of complementation.)

As mentioned in the introduction, we will consider the four symmetry classes
whose symmetry groups are subgroups of $\langle \rho, \kappa\tau \rangle$. 
In other words, CSPP's are invariant under $\rho$, TCPP's are invariant under
$\kappa\tau$, and CSTCPP's are invariant under both.

In addition to ordinary enumeration, we will also consider the weighted
enumeration of plane partitions, where the weight of a plane partition with
$n$ cubes is $q^n$.  Specifically, let $N_G(a,b,c)_q$ be the total weight of
all plane partitions in the given symmetry class.  The polynomials
$N(a,b,c)_q$, $N_\tau(a,a,b)_q$, and $N_\rho(a,a,a)_q$ have nice product
formulas, but we will only give a proof of the formula for $N(a,b,c)_q$.
(There is an alternate weighting whereby the weight of a plane partition is
$q^n$ if it has $n$ orbits; the corresponding polynomial $N'_\tau(a,a,b)_q$
is known to have a product formula and $N'_{\rho,\tau}(a,a,a)$ is conjectured
to do so, but we will not consider these interesting cases.)

\subsection{Product Formulas}

We define some products related to Stanley's formulas
\cite{Stanley:symmetries}.

If $a_1, a_2, \ldots, a_k$ are positive integers and $n$ is a
non-negative integer, define the \emph{ box product}
$C(a_1,a_2,\ldots,a_k;n)$ by:
$$C(a_1,a_2,\ldots,a_k;n) = \prod_{x_1 = 0}^{a_1-1}
\prod_{x_2 = 0}^{a_2-1} \ldots
\prod_{x_k = 0}^{a_k-1} \max(n-\sum_{i=1}^k x_i,1)$$
A typical box product is the product of the numbers shown in 
Figure~\ref{fboxprod}.  Define the \emph{ simplex product} $T(k,n)$ by:
$$T(k,n) = C(\overbrace{\infty,\infty,\ldots,\infty}^{k};n).$$ 
Thus, $T(1,n) = n!$, while 
$$T(2,n) = n!(n-1)!(n-2)!\ldots 2!1!.$$

\begin{fullfigure}{fboxprod}{Factors of the box product $C(6,4;7)$}
\psset{unit=.5cm}
\pspicture(-.5,-.5)(6.5,3.5)
\rput(0,0){7}\rput(0,1){6}\rput(0,2){5}\rput(0,3){4}
\rput(1,0){6}\rput(1,1){5}\rput(1,2){4}\rput(1,3){3}
\rput(2,0){5}\rput(2,1){4}\rput(2,2){3}\rput(2,3){2}
\rput(3,0){4}\rput(3,1){3}\rput(3,2){2}\rput(3,3){1}
\rput(4,0){3}\rput(4,1){2}\rput(4,2){1}\rput(4,3){1}
\rput(5,0){2}\rput(5,1){1}\rput(5,2){1}\rput(5,3){1}
\multips(-.5,-.5)(1,0){7}{\qline(0,0)(0,4)}
\multips(-.5,-.5)(0,1){5}{\qline(0,0)(6,0)}
\endpspicture
\end{fullfigure}

The simplex products can also be defined inductively by the following rules.
\begin{align}
T(0,n) &= n \nonumber \\
T(k,0) &= 1 \nonumber \\
T(k,n) &= T(k-1,n)T(k,n-1) \label{esimpind}
\end{align}
for positive integers $k$ and $n$, with $T(k,n) = 1$ for all other cases.

It follows immediately from the definition that
$$C(a_1,a_2,\ldots,a_k;n) = \frac{C(\infty,a_2,a_3,\ldots,a_k;n)}
    {C(\infty,a_2,a_3,\ldots,a_k;n-a_1)}.$$
Applying this formula to each variable in turn yields a multiplicative
inclusion-exclusion formula:
\begin{multline*}
C(a_1,a_2\ldots,a_k;n) = \\
    \begin{aligned}
    & T(k,n) T(k,n-a_1)^{-1}T(k,n-a_2)^{-1}\ldots T(k,n-a_k)^{-1} \\ 
    & T(k,n-a_1-a_2)T(k,n-a_1-a_3)\ldots \\
    & \quad\vdots \\
    & T(k,n-\sum_{i=1}^k a_i)^{(-1)^k}.
    \end{aligned}
\end{multline*}
In particular,
\begin{multline}
C(a,b,c;n) = \\
\frac{T(3,n) T(3,n-a-b) T(3,n-a-c) T(3,n-b-c)}{
T(3,n-a) T(3,n-b) T(3,n-c) T(3,n-a-b-c)}. \label{eincec}
\end{multline}

Similarly, define a box $q$-product by:
$$C(a_1,a_2,\ldots,a_k;n)_q = \prod_{0 \le x_1 \le a_i-1}
    \left(\max(n-\sum_{i=1}^k x_i,\ 1)\right)_q,$$
where $(n)_q = 1+q+\ldots+q^{n-1}$, and a simplex $q$-product by
$$T(k,n)_q = C(\infty,\ldots,\infty;n)_q.$$

The formula for $N(a,b,c)$ given by Stanley \cite{Stanley:symmetries}
and originally proved by MacMahon, is
$$N(a,b,c) = \frac{C(a,b,c;a+b+c-1)}{C(a,b,c;a+b+c-2)}.$$
Combining this equation with equation~(\ref{eincec}) and
equation~(\ref{esimpind}), we obtain the expression:
an expression 
$$ N(a,b,c) = 
    \frac{T(2,d-1)T(2,a-1)T(2,b-1)T(2,c-1)}
    {T(2,\chat-1)T(2,\bhat-1)T(2,\ahat-1)},$$
where, for abbreviation
\begin{align*}
\ahat &= b+c & \bhat &= a+c \\
\chat &= a+b & d &= a+b+c.
\end{align*}
This expression was suggested by Propp
\cite{Kuperberg:perdet} in the form
$$N(a,b,c) = \frac{H(a+b+c)H(a)H(b)H(c)}{H(a+b)H(a+c)H(b+c)},$$
where
$$H(n) = 1!2!\ldots(n-1)! = T(2,n-1)$$
is the hyperfactorial function.

For either formulation of MacMahon's enumeration, we can $q$-ify throughout.

\section{The enumeration}

\subsection{The permanent-determinant method}
\label{sperdet}

Following reference \cite{Kuperberg:perdet}, a plane partition in an $a
\times b \times c$ box is equivalent to a tiling of a certain hexagon
$H(a,b,c)$ by unit lozenges, where a \emph{ lozenge} is a rhombus with a
60 degree angle, as shown in Figure~\ref{ftiling}.  The hexagon
$H(a,b,c)$ has angles of 120 degrees and edge lengths of $a$, $b$, $c$,
$a$, $b$, and $c$, going clockwise around the perimeter.  The hexagon
has a unique tiling by unit equilateral triangles, and a unit lozenge
tiling is therefore equivalent to a perfect matching of a bipartite
graph $Z(a,b,c)$, shown in Figure~\ref{fgraph}, whose vertices are the
triangles and whose edges are given by adjacency.

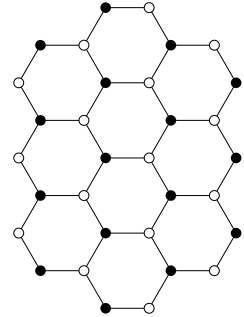
\begin{fullfigure}{fgraph}{The graph $Z(2,2,3)$}
\psset{xunit=.289cm,yunit=.5cm}
\pspicture(-6,0)(6,10)
\multips(-4,2)(0,2){4}{\qline(0,0)(2,0)}  
\multips(-1,1)(0,2){5}{\qline(0,0)(2,0)}
\multips( 2,2)(0,2){4}{\qline(0,0)(2,0)}
\multips(-4,2)(0,2){3}{\qline(0,0)(-1,1)} 
\multips(-1,1)(0,2){4}{\qline(0,0)(-1,1)}
\multips( 2,2)(0,2){4}{\qline(0,0)(-1,1)}
\multips( 5,3)(0,2){3}{\qline(0,0)(-1,1)}
\multips( 4,2)(0,2){3}{\qline(0,0)( 1,1)} 
\multips( 1,1)(0,2){4}{\qline(0,0)( 1,1)}
\multips(-2,2)(0,2){4}{\qline(0,0)( 1,1)}
\multips(-5,3)(0,2){3}{\qline(0,0)( 1,1)}
\multips(-5,3)(0,2){3}{\whitedot}\multips(-2,2)(0,2){4}{\whitedot}
\multips( 1,1)(0,2){5}{\whitedot}\multips( 4,2)(0,2){4}{\whitedot}
\multips( 5,3)(0,2){3}{\blackdot}\multips( 2,2)(0,2){4}{\blackdot}
\multips(-1,1)(0,2){5}{\blackdot}\multips(-4,2)(0,2){4}{\blackdot}
\endpspicture
\end{fullfigure}

The number of perfect matchings of an arbitrary bipartite graph is
given by the permanent of its bipartite adjacency matrix. Kasteleyn's
\emph{ permanent-determinant method} dictates that, if the graph is
planar, then there is a way to change the signs of the permanent to
convert it to a determinant.  (More accurately, Kasteleyn considered
the more general non-bipartite case and produced a Pfaffian; the
bipartite case was clarified by Percus \cite{Percus:dimer}.) Indeed,
the permanent of the adjacency matrix $A(a,b,c)$ of $Z(a,b,c)$ equals
the determinant up to a global sign, and the global sign is in any case
ambiguous because the rows and columns are unordered. One way to
demonstrate the equality is to note that two matchings differ by a
3-cycle (carried by a hexagon of $Z(a,b,c)$) if their plane partitions
differ by a cube.  Since a 3-cycle is an even permutation, it incurs no
sign change between the corresponding terms of $\Det A(a,b,c)$.  Since
any two plane partitions are connected by the operation of adding or
taking away an individual cube, it follows that all terms have the same
sign in the determinant.

More generally, suppose that $Z(a,b,c)$ is weighted with non-zero weights, and
that $M$ is its weighted adjacency matrix.  Then the non-vanishing terms in
the determinant of $M$ also correspond to matchings of $Z(a,b,c)$, but they
are no longer necessarily equal.  If, for every hexagon of $Z(a,b,c)$, the
weights $x_1,y_1,z_1,x_2,y_2$, and $z_2$ shown in Figure~\ref{fweights}
satisfy
$$x_1y_1z_1 = x_2y_2z_2,$$
then the non-vanishing terms of $\Det M$ are equal.  In this case, we say that
$M$ and the weighting are \emph{ Kasteleyn-flat}.  The number of matchings is
therefore 
\bqn (\Det M)/m,\label{eperdet}\eqn
where $m$ is the value of any single term.  This principle can be
extended to the $q$-enumeration problem:  Define the \emph{ Kasteleyn
curvature} of a hexagon as
$$\frac{x_1y_1z_1}{x_2y_2z_2}.$$
(This statistic is called curvature by analogy with the curvature of a
connection on a line bundle or the coboundary operation in homology; see
reference~\cite{Kuperberg:fun} for details.) If the Kasteleyn curvature is
$q$ everywhere and $t$ is the weight of the matching of the empty plane
partition, then
$$N(a,b,c)_q = (\Det M)/m.$$
In this case, the value of a term goes up by $q$ if we add a cube to the
corresponding plane partition.

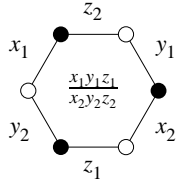
\begin{fullfigure}{fweights}{The weights that determine Kasteleyn curvature}
\psset{xunit=.433cm,yunit=.75cm}
\pspicture(-3,-1.5)(3,1.5)
\pcline(-1,-1)(-2, 0)\Aput{$y_2$}
\pcline(-2, 0)(-1, 1)\Aput{$x_1$}
\pcline(-1, 1)( 1, 1)\Aput{$z_2$}
\pcline( 1, 1)( 2, 0)\Aput{$y_1$}
\pcline( 2, 0)( 1,-1)\Aput{$x_2$}
\pcline( 1,-1)(-1,-1)\Aput{$z_1$}
\qdisk(-1,-1){3pt}\pscircle[fillstyle=solid,fillcolor=white](-2,0){3pt}
\qdisk(-1, 1){3pt}\pscircle[fillstyle=solid,fillcolor=white]( 1,1){3pt}
\qdisk( 2, 0){3pt}\pscircle[fillstyle=solid,fillcolor=white](1,-1){3pt}
\rput(0,0){$\frac{x_1y_1z_1}{x_2y_2z_2}$}
\endpspicture
\end{fullfigure}

The symmetry operations $\tau$, $\rho$, and $\kappa$ act on hexagons
$H(a,b,c)$ respectively as rotation by $120^\circ$, rotation by $180^\circ$,
and reflection about a diagonal. Therefore $\kappa\tau$ acts as reflection
about a bisector (a line which meets two opposite edges of the hexagon in
the middle), and in its action on $Z(a,a,2b)$, it fixes a row of
edges which separates the graph into two isomorphic subgraphs. (Note that
there are no $\kappa\tau$-invariant matchings of $Z(a,a,2b-1)$, or
equivalently $N_{\kappa\tau}(a,a,2b-1) = 0$.) The fixed edges must appear in
any invariant subgraph, so that TCPP's correspond to matchings of either
subgraph, denoted by $Z_{\kappa\tau}(a,a,2b)$.  Similarly, $\langle \rho,
\kappa\tau\rangle$ has three bisectors as lines of reflection, so that
CSTCPP's correspond to matchings of any of the six subgraphs
$Z_{\rho,\kappa\tau}(2a,2a,2a)$ that remain.  Since all the faces of these
graphs are also hexagons, the same analysis of Kasteleyn's method apply
fully, except that Kasteleyn curvature of $q$ is unnatural from the point of
view of plane partitions unless $q = \pm 1$.  We will only consider $q = 1$
for these cases in this paper.

A $\rho$-invariant matching of $Z(a,a,a)$ is simply a matching of the
quotient graph $Z_\rho(a,a,a) = Z(a,a,a)/\rho$, and again the same
analysis of Kasteleyn's method applies, except that the face in the
center is now a 2-gon.  This means that $Z_\rho(a,a,a)$ is not a simple
graph, the matrix $M$ is not uniquely determined by the weighting of
$Z_\rho(a,a,a)$, and the matchings are not bijective with the non-zero
terms in the determinant of $M$.  However, if the weight of each edge
is a separate variable, then the monomials in the expansion of $\det M$
are bijective with the matchings of $Z_\rho(a,a,a)$.  Taking this
approach, if the 2-gon has sides with weight $x_1$ and $x_2$, we define
its \emph{ Kasteleyn curvature} as $x_1/x_2$ and say that it is flat if
$x_1 = x_2$. As before, all terms in the determinant are equal if the
weighting is flat at every face.  More generally, the assignment of
Kasteleyn curvature which corresponds to the $q$-enumeration
$N_\rho(a,a,a)_q$ is $q^3$ for every hexagon and $q$ for the 2-gon.

\subsection{Representation theory of $\sl(2,\C)$ \label{ssl2}}

In this section, we will find a Kasteleyn-flat matrix $M$ for
$Z(a,b,c)$ in the representation theory of $\sl(2,\C)$.  We begin by
reviewing the finite-dimensional representation theory of $\sl(2,\C)$.

The Lie algebra $\sl(2,\C)$ of traceless $2 \times 2$ matrices has an
important basis $H$, $X$, and $Y$ given by the matrices:
\begin{align*}
H &= \begin{pmatrix} 1 & 0 \\ 0 & -1\end{pmatrix} \\
X &= \begin{pmatrix} 0 & 1 \\ 0 & 0\end{pmatrix} \\
Y &= \begin{pmatrix} 0 & 0 \\ 1 & 0\end{pmatrix}
\end{align*}
The entire Lie algebra structure can also be defined by linear
extension of the Lie bracket on this basis:
\begin{align*}
[H,X] &= 2X \\
[H,Y] &= -2Y \\
[X,Y] &= H
\end{align*}
Let $x$ and $y$ be the standard basis vector on which the above matrices act. 
For each positive integer $n$, $\sl(2,\C)$ has an irreducible
representation $V_n$ of dimension $n+1$ which can be described as the action
of $\sl(2,\C)$ on homogeneous polynomials in $x$ and $y$ of degree $n$. Let
$\alpha_n$ be the representation map from $\sl(2,\C)$ to $\End(V_n)$. 
To derive the action of $\sl(2,\C)$ in the 
the monomial basis $x^n, x^{n-1}y,x^{n-2}y^2,\ldots,y^n$,
it is convenient to express $X$, $Y$, and $H$ formally as:
\begin{align*}
X &= x \frac{\partial}{\partial y} \\
Y &= y \frac{\partial}{\partial x} \\
H &= x \frac{\partial}{\partial x} - y \frac{\partial}{\partial y}
\end{align*}
The corresponding matrices are then:
\begin{align*}
\alpha_n(H) &= \begin{pmatrix}
n & 0 & 0 & \cdots & 0 \\
0 & n-2 & 0 & \cdots & 0 \\
0 & 0 & n-4 & \cdots & 0 \\
\vdots & \vdots & \vdots & \ddots & \vdots \\
0 & 0 & 0 & \cdots & -n \end{pmatrix} \\
\alpha_n(X) &= \begin{pmatrix}
0 & 1 & 0 & 0 & \cdots & 0 \\
0 & 0 & 2 & 0 & \cdots & 0 \\
0 & 0 & 0 & 3 & \cdots & 0 \\
\vdots & \vdots & \vdots & \vdots & \ddots & \vdots \\
0 & 0 & 0 & 0 & \cdots & n \\
0 & 0 & 0 & 0 & \cdots & 0 \end{pmatrix} \\
\alpha_n(Y) &= \begin{pmatrix}
0 & 0 & 0 & \cdots & 0 & 0 \\
n & 0 & 0 & \cdots & 0 & 0 \\
0 & n-1 & 0 & \cdots & 0 & 0 \\
0 & 0 & n-2 & \cdots & 0 & 0 \\
\vdots & \vdots & \vdots & \ddots & \vdots & \vdots \\
0 & 0 & 0 & \cdots & 1 & 0 \end{pmatrix}
\end{align*}

The monomial basis is also known as the (dual) \emph{ weight basis}, and
we rename the basis vectors $e_n,e_{n-2},\ldots,e_{-n}$, where the
eigenvalue of $e_i$ with respect to $H$ is $i$.

\begin{theorem} (Clebsch-Gordan)  Any finite-dimensional, irreducible
representation of $\sl(2,\C)$ is isomorphic to a direct sum of $V_n$'s.
In particular, their tensor products decompose according to the
equation:
$$V_n \tensor V_k \cong V_{n+k} \oplus V_{n+k-2} \oplus V_{n+k-4} \oplus
\ldots \oplus V_{|n-k|}.$$
\end{theorem}

We assume the first half of the theorem without proof.  To understand
the Clebsch-Gordan decomposition of a tensor product, we 
recall that if $L$ is an element of an
arbitrary Lie algebra, its action on the tensor product of two
representations $\alpha$ and $\beta$ is given by:
\begin{equation}
(\alpha \tensor\beta)(L) = \alpha(L) \tensor I + I \tensor \beta(L),
\label{etensor}
\end{equation}
where $I$ is the identity matrix.

Given the first half of the theorem, the following is an argument for
the Clebsch-Gordan formula:  Let $S(H)$ and $S(\sl(2,\C))$ be the
semi-ring of reducible representations of $H$ with integer eigenvalues
and the semi-ring of representations of $\sl(2,\C)$, respectively.  The
elements of these semi-rings are representations and the algebraic
structure is given by $\oplus$ and $\otimes$.   The semi-rings extend
to \emph{ Grothendieck rings} $R(H)$ and $R(\sl(2,\C))$ by introducing
subtraction. Since each $V_n$ restricts to a reducible representation
of $H$, $S(\sl(2,\C)) \subset S(H)$ and $R(\sl(2,\C)) \subset R(H)$. The ring
$R(H)$, as an abelian group, is generated by one-dimensional
representations $E_n$ with eigenvalue $n$ for some integer $n$.  By
equation~(\ref{etensor}), eigenvalues add under tensor product, and
therefore there exists an isomorphism $\ch:R(H) \to \Z[t,1/t]$ given by
the formula $\ch(E_n) = t^n$. This isomorphism has a restriction
$\ch:R(\sl(2,\C)) \to \Z[t,1/t]$ that we will call the \emph{ character
map}. By the form of $\alpha_n(H)$,
\bqn \ch(V_n) = t^n + t^{n-2} + \ldots + t^{-n}. \label{ewt} \eqn

Because the $\ch(V_n)$'s are linearly independent, $\ch$ is injective. The
Clebsch-Gordan formula then follows easily by applying $\ch$ to both sides. 
Indeed, the assertion that the character map is injective on
$S(\sl(2,\C))$ is a convenient restatement of the Clebsch-Gordan theorem.

Since the matrices for $\alpha_n(H)$, $\alpha_n(X)$, and $\alpha_n(Y)$ are
nearly permutation matrices, they can be described graphically as in
Figure~\ref{fsltwoirred}.  Each dot in the figure represents a basis vector
which is an eigenvector of $\alpha_n(H)$ with the given eigenvalue.  The maps
$\alpha_n(X)$ and $\alpha_n(Y)$ send each basis vector to another basis
vector times the given scalar factor, except for the vector at the right
(resp. left), which is in the kernel of $\alpha_n(X)$ (resp. $\alpha_n(Y)$).

\begin{fullfigure}{fsltwoirred}{An irreducible representation of $\sl(2)$}
\psset{unit=1cm,nodesep=6pt}
{\small\pspicture(-4.5,-2)(3.5,2)
\multips(-3,0)(1,0){7}{\qdisk(0,0){2pt}}
\pcarc[arcangle=25]{->}(-3,0)(-2,0)\Aput{6}
\pcarc[arcangle=25]{->}(-2,0)(-1,0)\Aput{5}
\pcarc[arcangle=25]{->}(-1,0)( 0,0)\Aput{4}
\pcarc[arcangle=25]{->}( 0,0)( 1,0)\Aput{3}
\pcarc[arcangle=25]{->}( 1,0)( 2,0)\Aput{2}
\pcarc[arcangle=25]{->}( 2,0)( 3,0)\Aput{1}
\pcarc[arcangle=25]{->}( 3,0)( 2,0)\Aput{6}
\pcarc[arcangle=25]{->}( 2,0)( 1,0)\Aput{5}
\pcarc[arcangle=25]{->}( 1,0)( 0,0)\Aput{4}
\pcarc[arcangle=25]{->}( 0,0)(-1,0)\Aput{3}
\pcarc[arcangle=25]{->}(-1,0)(-2,0)\Aput{2}
\pcarc[arcangle=25]{->}(-2,0)(-3,0)\Aput{1}
\pcline[linestyle=dashed](-3,0)(-3,-1.2)\rput[t](-3,-1.2){-6}
\pcline[linestyle=dashed](-2,0)(-2,-1.2)\rput[t](-2,-1.2){-4}
\pcline[linestyle=dashed](-1,0)(-1,-1.2)\rput[t](-1,-1.2){-2}
\pcline[linestyle=dashed]( 0,0)( 0,-1.2)\rput[t]( 0,-1.2){0}
\pcline[linestyle=dashed]( 1,0)( 1,-1.2)\rput[t]( 1,-1.2){2}
\pcline[linestyle=dashed]( 2,0)( 2,-1.2)\rput[t]( 2,-1.2){4}
\pcline[linestyle=dashed]( 3,0)( 3,-1.2)\rput[t]( 3,-1.2){6}
{\normalsize\rput[t](-4,-1.2){$H$}\rput(-4,-.4){$Y$}\rput(-4,.4){$X$}}
\endpspicture}
\end{fullfigure}

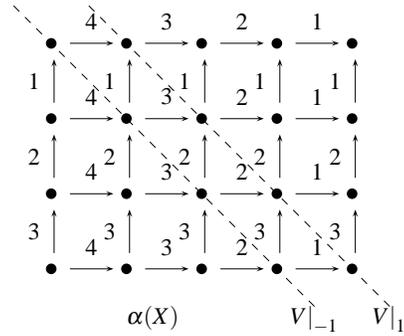
\begin{fullfigure}{fsltwontensork}{The representation $V_4 \tensor V_3$}
\psset{unit=1cm,nodesep=6pt}
\pspicture(-.5,-1.5)(5.5,5)
{\small
\multips(0,0)(1,0){5}{
\qdisk(0,0){2pt}\qdisk(0,1){2pt}\qdisk(0,2){2pt}\qdisk(0,3){2pt}}
\multirput(0,0)(0,1){4}{
\pcline{->}(0,0)(1,0)\Aput{4}\pcline{->}(1,0)(2,0)\Aput{3}
\pcline{->}(2,0)(3,0)\Aput{2}\pcline{->}(3,0)(4,0)\Aput{1}}
\multirput(0,0)(1,0){5}{
\pcline{->}(0,0)(0,1)\Aput{3}\pcline{->}(0,1)(0,2)\Aput{2}
\pcline{->}(0,2)(0,3)\Aput{1}}}
\psline[linestyle=dashed](-.5,3.5)(3.5,-.5)
\psline[linestyle=dashed](.5,3.5,)(4.5,-.5)
\rput[t](3.5,-.5){$V|_{-1}$}
\rput[t](4.5,-.5){$V|_1$}
\rput[tl](1,-.5){$\alpha(X)$}
\endpspicture
\end{fullfigure}

As a warm-up for the main construction, consider the action of $X$,
$Y$, and $H$ on a tensor product representation $V = V_n \tensor V_k$
in the tensor product basis  $\{e_1 \tensor e_j\}$.  If $\alpha$ is the
representation map, the actions of $\alpha(H)$ and $\alpha(X)$ are
represented diagrammatically by Figure~\ref{fsltwontensork}.  Each dot
is a basis vector as before, and diagonal strings span eigenspaces of
$\alpha(H)$ with the given eigenvectors, because $e_i \tensor e_j$ is
an eigenvector of $\alpha(H)$ with eigenvalue $i+j$. By
equation~(\ref{etensor}), most columns of the matrix for $\alpha(X)$
have two terms, which are given by the arrows. (Technically speaking,
one might say that the diagram is a directed, weighted graph and the
matrix for $\alpha(X)$ is the asymmetric weighted adjacency matrix of
the graph.)  Let $V|_\lambda$ be the eigenspace of $\alpha(H)$ with
eigenvalue $\lambda$.  Then the matrix for $\alpha(X)$ can be divided
into blocks $\alpha(X)|_\lambda:V|_\lambda \to V|_{\lambda+2}$ which
are maps between adjacent eigenspaces.  For example, in
Figure~\ref{fsltwontensork}, $\alpha(X)|_{-1}$, whose domain and target
are delineated, has the matrix:

$$\begin{pmatrix}
4 & 1 & 0 & 0 \\
0 & 3 & 2 & 0 \\
0 & 0 & 2 & 3 \\
0 & 0 & 0 & 1 \end{pmatrix}$$

Consider the graph $Z(a,b,c)$ as in 
Section~\ref{sperdet}, let 
$$V = V_{\chat-1} \tensor V_{\bhat-1} \tensor V_{\ahat-1},$$
and let $\alpha$ be the representation map.  The tensor product basis  $\{e_i
\tensor e_j \tensor e_k\}$ for this representation can be depicted as a
3-dimensional rectangular box of points. Unfortunately, it is hard to show all of
the points in a 2-dimensional picture. However, $\alpha(H)$ is diagonal as
before, and the basis vectors lying in each eigenspace constitute a 2-dimensional
slice of the box, because each $e_i \tensor e_j \tensor e_k$ is an eigenvector
with eigenvalue $i+j+k$.  These slices are given as a sequence for $V_4 \tensor
V_4 \tensor V_5$ in Figure~\ref{fsltwoslices}, together with the eigenvalue of
$\alpha(H)$ for each slice. We will consider the slices $V|_{-1}$
and $V|_1$, which are given separately in Figure~\ref{fsltwodots},
together with points corresponding to their basis vectors.

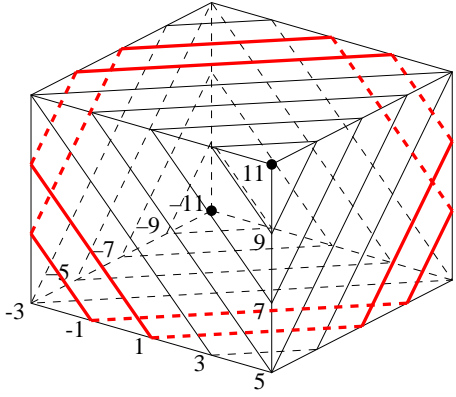
\begin{fullfigure}{fsltwoslices}{The representation $V_4 \tensor V_4 \tensor V_5$}
\SpecialCoor
\psset{unit=1cm}
{\small \pspicture(-4,-3)(4,3)
\pspolygon(\td{4}{0}{0})(\td{4}{0}{3})(\td{0}{0}{3})(\td{0}{4}{3})%
(\td{0}{4}{0})(\td{4}{4}{0})
\psline(\td{4}{4}{0})(\td{4}{4}{3})
\psline(\td{4}{0}{3})(\td{4}{4}{3})
\psline(\td{0}{4}{3})(\td{4}{4}{3})
\psline[linestyle=dashed](\td{0}{0}{0})(\td{0}{0}{3})
\psline[linestyle=dashed](\td{0}{0}{0})(\td{0}{4}{0})
\psline[linestyle=dashed](\td{0}{0}{0})(\td{4}{0}{0})
\qdisk(\td{4}{4}{3}){2pt}
\qdisk(\td{0}{0}{0}){2pt}
\rput[br](\td{0}{0}{0}){--11\ }
\rput[br](\td{1}{0}{0}){--9\ }
\rput[br](\td{2}{0}{0}){--7\ }
\rput[br](\td{3}{0}{0}){--5\ }
\pspolygon(\td{3}{4}{3})(\td{4}{3}{3})(\td{4}{4}{2})
\pspolygon(\td{2}{4}{3})(\td{4}{2}{3})(\td{4}{4}{1})
\pspolygon(\td{1}{4}{3})(\td{4}{1}{3})(\td{4}{4}{0})
\psline(\td{3}{4}{0})(\td{0}{4}{3})(\td{4}{0}{3})(\td{4}{3}{0})
\psline[linestyle=dashed](\td{3}{4}{0})(\td{4}{3}{0})
\psset{linestyle=dashed}
\pspolygon(\td{1}{0}{0})(\td{0}{1}{0})(\td{0}{0}{1})
\pspolygon(\td{2}{0}{0})(\td{0}{2}{0})(\td{0}{0}{2})
\pspolygon(\td{3}{0}{0})(\td{0}{3}{0})(\td{0}{0}{3})
\psline(\td{1}{0}{3})(\td{4}{0}{0})(\td{0}{4}{0})(\td{0}{1}{3})
\psline[linestyle=solid](\td{0}{1}{3})(\td{1}{0}{3})
\psset{linewidth=1.25pt,linecolor=red}
\psline[linestyle=dashed](\td{2}{0}{3})(\td{4}{0}{1})
 \psline[linestyle=solid](\td{4}{0}{1})(\td{4}{1}{0})
\psline[linestyle=dashed](\td{4}{1}{0})(\td{1}{4}{0})
 \psline[linestyle=solid](\td{1}{4}{0})(\td{0}{4}{1})
\psline[linestyle=dashed](\td{0}{4}{1})(\td{0}{2}{3})
 \psline[linestyle=solid](\td{0}{2}{3})(\td{2}{0}{3})
\psline[linestyle=dashed](\td{3}{0}{3})(\td{4}{0}{2})
 \psline[linestyle=solid](\td{4}{0}{2})(\td{4}{2}{0})
\psline[linestyle=dashed](\td{4}{2}{0})(\td{2}{4}{0})
 \psline[linestyle=solid](\td{2}{4}{0})(\td{0}{4}{2})
\psline[linestyle=dashed](\td{0}{4}{2})(\td{0}{3}{3})
 \psline[linestyle=solid](\td{0}{3}{3})(\td{3}{0}{3})
\psset{linestyle=solid,linewidth=.25pt,linecolor=black}
\rput[tr](\td{4}{0}{0}){-3\ }
\rput[tr](\td{4}{1}{0}){-1\ }
\rput[tr](\td{4}{2}{0}){1\ }
\rput[tr](\td{4}{3}{0}){3\ }
\rput[tr](\td{4}{4}{0}){5\ }
\rput[tr](\td{4}{4}{1}){7\ }
\rput[tr](\td{4}{4}{2}){9\ }
\rput[tr](\td{4}{4}{3}){11\ }

\endpspicture}
\NormalCoor
\end{fullfigure}

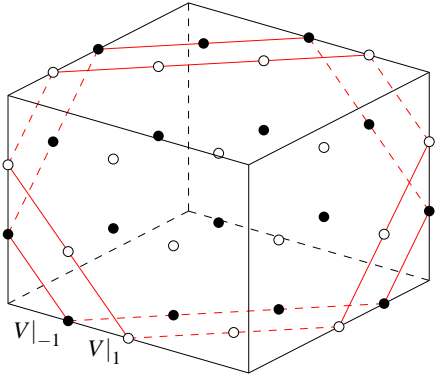
\begin{fullfigure}{fsltwodots}{Basis vectors in $V|_{-1}$ and $V|_1$}
\SpecialCoor
\psset{unit=1cm}
{\small \pspicture(-4,-3)(4,3)
\pspolygon(\td{4}{0}{0})(\td{4}{0}{3})(\td{0}{0}{3})(\td{0}{4}{3})%
(\td{0}{4}{0})(\td{4}{4}{0})
\psline[linestyle=dashed](\td{0}{0}{0})(\td{0}{0}{3})
\psline[linestyle=dashed](\td{0}{0}{0})(\td{0}{4}{0})
\psline[linestyle=dashed](\td{0}{0}{0})(\td{4}{0}{0})
\psset{linecolor=red}
\psline[linestyle=dashed](\td{2}{0}{3})(\td{4}{0}{1})
 \psline[linestyle=solid](\td{4}{0}{1})(\td{4}{1}{0})
\psline[linestyle=dashed](\td{4}{1}{0})(\td{1}{4}{0})
 \psline[linestyle=solid](\td{1}{4}{0})(\td{0}{4}{1})
\psline[linestyle=dashed](\td{0}{4}{1})(\td{0}{2}{3})
 \psline[linestyle=solid](\td{0}{2}{3})(\td{2}{0}{3})
\psline[linestyle=dashed](\td{3}{0}{3})(\td{4}{0}{2})
 \psline[linestyle=solid](\td{4}{0}{2})(\td{4}{2}{0})
\psline[linestyle=dashed](\td{4}{2}{0})(\td{2}{4}{0})
 \psline[linestyle=solid](\td{2}{4}{0})(\td{0}{4}{2})
\psline[linestyle=dashed](\td{0}{4}{2})(\td{0}{3}{3})
 \psline[linestyle=solid](\td{0}{3}{3})(\td{3}{0}{3})
\psset{linecolor=black}
\multips(\td{2}{0}{3})(\td{-1}{1}{0}){3}{\blackdot}
\multips(\td{3}{0}{2})(\td{-1}{1}{0}){4}{\blackdot}
\multips(\td{4}{0}{1})(\td{-1}{1}{0}){5}{\blackdot}
\multips(\td{4}{1}{0})(\td{-1}{1}{0}){4}{\blackdot}
\multips(\td{2}{4}{0})(\td{1}{-1}{0}){3}{\whitedot}
\multips(\td{1}{4}{1})(\td{1}{-1}{0}){4}{\whitedot}
\multips(\td{0}{4}{2})(\td{1}{-1}{0}){5}{\whitedot}
\multips(\td{0}{3}{3})(\td{1}{-1}{0}){4}{\whitedot}
\psline(\td{4}{4}{0})(\td{4}{4}{3})
\psline(\td{4}{0}{3})(\td{4}{4}{3})
\psline(\td{0}{4}{3})(\td{4}{4}{3})

\rput[tr](\td{4}{1}{0}){$V|_{-1}$\ }
\rput[tr](\td{4}{2}{0}){$V|_1$\ }
\endpspicture}
\NormalCoor
\end{fullfigure}

Most columns of the matrix for $\alpha(X)$ now have three terms, and the
matrix divides into submatrices $\alpha(X)|_\lambda:V|_\lambda \to
V|_{\lambda+2}$.  Consider the map $\alpha(X)|_{-1}$, which is diagrammed in
Figure~\ref{fxminone}.  If we replace the arrows in this diagram by
unoriented edges, we see that $\alpha(X)|_{-1}$ is a weighted, bipartite
adjacency matrix of the graph $Z(a,b,c)$.  Moreover, the weighting is
Kasteleyn-flat, because the weights of opposite edges of each face are equal.
We have proved the following proposition:

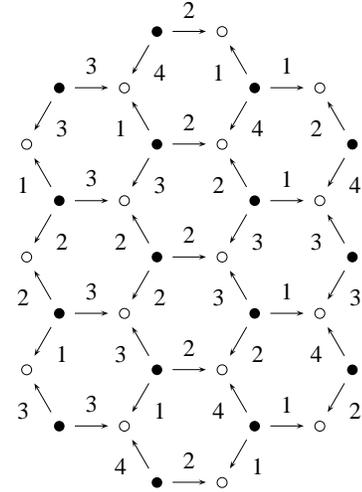
\begin{fullfigure}{fxminone}{The weights on $Z(2,2,3)$ coming from $\alpha(X)|_{-1}$}
\psset{xunit=.433cm,yunit=.75cm,nodesep=6pt}
\small{\pspicture(-6,0)(6,10)
\multirput(-4,2)( 0,2){4}{\pcline{->}(0,0)( 2,0)\Aput{3}} 
\multirput(-1,1)( 0,2){5}{\pcline{->}(0,0)( 2,0)\Aput{2}}
\multirput( 2,2)( 0,2){4}{\pcline{->}(0,0)( 2,0)\Aput{1}}
\multirput(-1,1)( 3,1){3}{\pcline{->}(0,0)(-1,1)\Aput{4}} 
\multirput(-4,2)( 3,1){4}{\pcline{->}(0,0)(-1,1)\Aput{3}}
\multirput(-4,4)( 3,1){4}{\pcline{->}(0,0)(-1,1)\Aput{2}}
\multirput(-4,6)( 3,1){3}{\pcline{->}(0,0)(-1,1)\Aput{1}}
\multirput( 1,1)(-3,1){3}{\pcline{<-}(0,0)( 1,1)\Bput{1}} 
\multirput( 4,2)(-3,1){4}{\pcline{<-}(0,0)( 1,1)\Bput{2}}
\multirput( 4,4)(-3,1){4}{\pcline{<-}(0,0)( 1,1)\Bput{3}}
\multirput( 4,6)(-3,1){3}{\pcline{<-}(0,0)( 1,1)\Bput{4}}
\multips(-5,3)(0,2){3}{\whitedot}\multips(-2,2)(0,2){4}{\whitedot}
\multips( 1,1)(0,2){5}{\whitedot}\multips( 4,2)(0,2){4}{\whitedot}
\multips( 5,3)(0,2){3}{\blackdot}\multips( 2,2)(0,2){4}{\blackdot}
\multips(-1,1)(0,2){5}{\blackdot}\multips(-4,2)(0,2){4}{\blackdot}
\endpspicture}
\end{fullfigure}

\begin{proposition} If the representation $V_{\hat-1} \tensor V_{\bhat-1} \tensor
V_{\chat-1}$ of $\sl(2,\C)$ has representation map $\alpha$ and is given with the
basis which is the tensor product of the weight bases of the three factors,
then the matrix for $\alpha(X)|_{-1}$, the restriction of $\alpha(X)$ to the
eigenspace of $\alpha(H)$ with eigenvalue $-1$, is a Kasteleyn-flat adjacency
matrix for $Z(a,b,c)$.
\end{proposition}

Applying equation~(\ref{eperdet}),
\bqn N(a,b,c) = (\det(\alpha(X)|_{-1}))/m,\label{enformula}\eqn
where $m$ is the value of any single term in the expansion of the
determinant.  Indeed, in this case all terms are not only equal, but
are products of the same weights.  The product of the weights of 
the horizontal edges in any matching is:
$$\prod_{i=1}^a \prod_{j=1}^b (i+j-1) = C(a,b;a+b-1)$$
Therefore
\bqn t = \pm C(a,b;\chat-1)C(a,c;\bhat-1)C(b,c;\ahat-1). \label{etval} \eqn

The map $\alpha(X)|_{-1}$ is a linear transformation between different
vector spaces.  For the purpose of computing determinants, it is
easier to work with an endomorphism of a single vector space.  Observe
that the restriction $\alpha(Y)|_1:V|_1 \to V|_{-1}$ has all of the
same properties as $\alpha(X)|_{-1}$, and that the composition
$\alpha(Y)|_1\alpha(X)|_{-1}$ is an endomorphism of $V|_{-1}$. In
particular,
$$N(a,b,c)^2 = \det(\alpha(Y)\alpha(X)|_{-1})/m^2.$$

Since $\alpha(Y)\alpha(X)|_{-1}$ comes from the action of $\sl(2,\C)$
on $V$, it acts on each summand of a direct-sum decomposition of $V$
separately.   Therefore such a decomposition diagonalizes
$\alpha(Y)\alpha(X)|_{-1}$ and can be used to find its determinant.

\begin{lemma} The map 
$$\alpha_{2n-1}(Y)\alpha_{2n-1}(X)|_{-1}:V_{2n-1}|_{-1} \to V_{2n-1}|_{-1}$$
is a map from a 1-dimensional vector space to itself and its effect
is multiplication by $n^2$.  \label{lnsquared}
\end{lemma}

Compare with Figure~\ref{fsltwoirred} or
with the definition of $\alpha_n(X)$ and $\alpha_n(Y)$.

In general if $\alpha$ and $\beta$ are any two representation of $\sl(2,\C)$,
\begin{multline*}
\det((\alpha+\beta)(Y)(\alpha+\beta)(X)|_{-1}) = \\
    \det(\alpha(Y)\alpha(X)|_{-1})\det(\beta(Y)\beta(X)|_{-1}).
\end{multline*}
It follows that the assignment
$$\alpha \mapsto \det(\alpha(Y)\alpha(X)|_{-1})$$
defines an abelian-group homomorphism from the additive subgroup of the
Grothendieck ring $R(\sl(2,\C))$ spanned by $V_{2n-1}$'s to the
multiplicative group $\Q^*$ of non-zero rationals.  This homomorphism
factors through the character map $\ch$ to yield a map $D$ such that
$$D(\ch(V)) = \det(\alpha(Y)\alpha(X))|_{-1})$$
for any odd-weight representation $V$.  By equation~(\ref{ewt}), we
can define $D$ by
$$D(t^{2n-1}) = \frac{n^2}{(n-1)^2}$$
for $n>1$ and $D(t^{2n-1}) = 1$ otherwise.

Our particular $V$ has character
$$\ch(V) = \sum_{(i,j,k) \in B} t^{2(i+j+k)-2d-3},$$
where
$$B = \{1,\ldots, \chat\} \times \{1,\ldots,\bhat\} \times \{1,\ldots,\ahat\}.$$
The map $D$ therefore yields
\begin{align*}
\det(\alpha(Y)\alpha(X)|_{-1}) &= \prod_{(i,j,k) \in B'}
    \frac{(i+j+k-d-1)^2}{(i+j+k-d-2)^2} \\
    &= \frac{C(\chat,\bhat,\ahat;d-1)^2}{C(\chat,\bhat,\ahat;d-2)^2},
\end{align*}
where $B' \subset B$ is the set of triples $(i,j,k)$ such
that
$$i+j+k \ge d+3.$$
Therefore
$$\det(\alpha(X)|_{-1}) =
\pm\frac{C(\chat,\bhat,\ahat;d-1)}{C(\chat,\bhat,\ahat;d-2)}$$
in the given basis. Combining this equation with equations~(\ref{etval})
and (\ref{enformula}), we obtain
$$N(a,b,c) = \frac{C(\chat,\bhat,\ahat;d-1)}
    {C(\chat,\bhat,\ahat;d-2) \prod C(a,b;\chat-1)}.$$
Here, the symbol ``$\prod$'' denotes repetition of the succeeding
factor with cyclic permutation of $a$, $b$, and $c$.

We apply the inclusion-exclusion formula for each cube product
and the inductive formula for the resulting simplex products:
\begin{align*}
N(a,b,c) &= \frac{T(3,d-1)\prod T(3,a-2)\prod T(2,a-1)^2}
    {T(3,d-2)\prod T(3,a-1) \prod T(2,\chat-1)} \\
    &= \frac{T(2,d-1)\prod T(2,a-1)}{\prod T(2,\chat-1)}.
\end{align*}
This is Propp's expression for $N(a,b,c)$.

\begin{remark}
Observe that $\alpha(X)|_{-1}$ is a diagonalizable,
integer matrix with integer eigenvalues.  The author considered a
direct analysis of this matrix without using any representation theory
other than the knowledge that the matrix can be decomposed using
rational linear algebra. Unfortunately, no obvious pattern for the
eigenvectors appeared in small examples.  One way to express the
coefficients of the eigenvectors with the aid of representation theory is
to use Racah $3j$-symbols, which are complicated and have no obvious
derivation using basic linear algebra.
\end{remark}

\section{Quantum representation theory of $\sl(2,\C)$}

The ordinary representation theory of $\sl(2,\C)$ admits a $q$-analogue,
the \emph{ quantum representation theory}, which can be used to
$q$-enumerate plane partitions.  A proper algebraic treatment of the
quantum representation theory would present it as the representation
theory of a Hopf algebra $U_q(\sl(2,\C))$ \cite{Drinfeld:quantum}.
Briefly, an abstract representation category, or the representation
theory of a Hopf algebra, is a category of vector spaces on which some
associative algebra or some matrices act, together with suitable
definitions of trivial representations, tensor-product representation,
and dual representations.  The most apparent difference between such a
category and the ordinary representation theory of a Lie algebra is
that $A \tensor B$ and $B \tensor A$ need not be isomorphic, and if
they are isomorphic, they need not be canonically isomorphic.  For
brevity, we will avoid the definition of $U_q(\sl(2,\C))$ and give a
direct, computational definition of quantum representations of
$\sl(2,\C)$.

Let $h$ be a complex number, and define $q = \exp(h)$.  More generally,
for any number or matrix $A$,  define $q^A =\exp(hA)$.  Define the
\emph{ bracket} of a number or matrix by
$$[A] = \frac{q^{A/2} - q^{-A/2}}{q^{1/2} - q^{-1/2}}.$$
The Laurent polynomial $[n]$ is called a \emph{ quantum integer} and is
related to the $q$-integer $(n)_q$ by
$$[n] = q^{(1-n)/2} (n)_q.$$

A \emph{quantum representation} of $\sl(2,\C)$ is defined to be a vector space
$V$ and a function
$\alpha:\{H,X,Y\} \to \End(V)$ such that
\begin{align*}
[\alpha(H),\alpha(X)] &= 2\alpha(X), \\
[\alpha(H),\alpha(Y)] &= -2\alpha(Y), \\
[\alpha(X),\alpha(Y)] &= [\alpha(H)].
\end{align*}
The tensor product of two representation maps $\alpha$ and $\beta$ is defined
by the equations
\begin{align*}
(\alpha\tensor\beta)(X) &= \alpha(X)\tensor q^{\beta(H)/4}
          + q^{-\alpha(H)/4}\tensor\beta(X), \\
(\alpha\tensor\beta)(Y) &= \alpha(Y)\tensor q^{\beta(H)/4}
          + q^{-\alpha(H)/4}\tensor\beta(Y), \\
(\alpha\tensor\beta)(H) &= \alpha(H)\tensor I + I \tensor \beta(H).
\end{align*}
These equations are analyzed by Drinfel'd \cite{Drinfeld:quantum}, but 
the diligent reader can check directly that the tensor product
matrices again form a valid quantum representation.  As a first
step, $[(\alpha\tensor\beta)(H)]$ involves $q^{(\alpha \tensor \beta)(H)/2}$,
and
\begin{align*}
q^{(\alpha\tensor\beta)(H)/2}
    &= q^{\alpha(H)/2 \tensor I + I \tensor \beta(H)/2} \\
    &= q^{\alpha(H)/2 \tensor I}q^{I \tensor \beta(H)/2} \\
    &= q^{\alpha(H)/2} \tensor q^{\beta(H)/2},
\end{align*}
since the exponential of a sum of commuting matrices is the product of
the exponentials.

The classical representation $V_n$ deforms to a quantum representation
by the formulas
\begin{align*}
\alpha_n(H) &= \begin{pmatrix}
    n & 0 & 0 & \cdots & 0 \\
    0 & n-2 & 0 & \cdots & 0 \\
    0 & 0 & n-4 & \cdots & 0 \\
    \vdots & \vdots & \vdots & \ddots & \vdots \\
0 & 0 & 0 & \cdots & -n \end{pmatrix}, \\
\alpha_n(X) &= \begin{pmatrix}
    0 & [1] & 0 & 0 & \cdots & 0 \\
    0 & 0 & [2] & 0 & \cdots & 0 \\
    0 & 0 & 0 & [3] & \cdots & 0 \\
    \vdots & \vdots & \vdots & \vdots & \ddots & \vdots \\
    0 & 0 & 0 & 0 & \cdots & [n] \\
    0 & 0 & 0 & 0 & \cdots & 0 \end{pmatrix}, \\
\alpha_n(Y) &= \begin{pmatrix}
    0 & 0 & 0 & \cdots & 0 & 0 \\[0pt]
    [n] & 0 & 0 & \cdots & 0 & 0 \\
    0 & [n-1] & 0 & \cdots & 0 & 0 \\
    0 & 0 & [n-2] & \cdots & 0 & 0 \\
    \vdots & \vdots & \vdots & \ddots & \vdots & \vdots \\
    0 & 0 & 0 & \cdots & [1] & 0 \end{pmatrix}.
\end{align*}
Finally, although in the quantum representation theory there is no canonical
isomorphism between $V_n \tensor V_k$ and $V_k \tensor V_n$, they are
isomorphic.  Indeed, the Clebsch-Gordan theorem generalizes to the quantum
representation theory.

As before, let
$$V = V_{\chat-1} \tensor V_{\bhat-1} \tensor V_{\ahat-1},$$
and define $V|_\lambda$ to be the eigenspace of $\alpha(H)$ with eigenvalue
$\lambda$ in this representation.  In the tensor product basis, the matrix for
$\alpha(X)|_{-1}$ is again a weighted adjacency matrix for
$Z(a,b,c)$.  In the quantum case, the weighting  has Kasteleyn curvature $q$ at
every hexagonal face.  Therefore
$$N(a,b,c)_q = (\det \alpha(X)|_{-1})/m_q,$$
where $m_q$ is the weight of the matching corresponding to the empty
plane partition.

The rest of the derivation is a $q$-analogue (or quantization)
of Section~\ref{ssl2}.  The quantity $m_q$ is given by
$$m_q = \pm q^{-abc/2} C(a,b;\chat-1)_qC(a,c;\bhat-1)_qC(b,c;\ahat-1)_q,$$
and the Clebsch-Gordan decomposition of $V$ yields
$$\det(\alpha(X)|_{-1})= \pm q^{-abc/2}
    \frac{C(\chat,\bhat,\ahat;d-1)_q}{C(\chat,\bhat,\ahat;d-2)_q}.$$
We can $q$-ify the derivation at the end of Section~\ref{ssl2}
to prove the $q$-analogue of equation~(\ref{eincec}):
$$
N(a,b,c)_q = \frac{T(2,d-1)_qT(2,a-1)_qT(2,b-1)_qT(2,c-1)_q}
    {T(2,\chat-1)_qT(2,\bhat-1)_qT(2,\ahat-1)_q}.$$

\section{Other determinantal symmetry classes \label{sothersl2}}

The method of Section~\ref{ssl2} naturally extend to the enumeration
of three other symmetry classes of plane partitions.

\begin{theorem} A Kasteleyn-flat, weighted adjacency matrix for each of the
graphs $Z_{\kappa\tau}(a,a,2b)$, $Z_\rho(a,a,a)$, and
$Z_{\rho,\kappa\tau}(2a,2a,2a)$ arises as $\alpha_{\kappa\tau}(X)|_{-1}$,
$\alpha_{\rho}(X)|_{-1}$, and $\alpha_{\rho,\kappa\tau}(X)|_{-1}$ in
representations $\alpha_{\kappa\tau}$, $\alpha_{\rho}$,
and $\alpha_{\rho,\kappa\tau}$ of $\sl(2,\C)$, with bases formed from weight
bases of irreducible representations.
\end{theorem}

Recall the weighting of $Z(a,a,2b)$ induced by $\alpha(X)|_{-1}$ acting on
$V = V_{a+2b-1}\tensor V_{a+2b-1}\tensor V_{2a-1}$.  The symmetry $\kappa\tau$
preserves this weighting, and is induced by the linear transformation
$$L_{\kappa\tau}\in\End(V_{a+2b-1}\tensor V_{a+2b-1}\tensor V_{2a-1})$$
given by
$$L_{\kappa\tau}(e_i \tensor e_j \tensor e_k) = e_j \tensor e_i \tensor e_k.$$
Define $V_{\kappa\tau}$ to be the eigenspace of $L_{\kappa\tau}$ with
eigenvalue $-1$, and choose vectors of the form $e_i \tensor e_j \tensor e_k
- e_j \tensor e_i \tensor e_k$ with $i<j$ as a basis.  Since $L_{\kappa\tau}$
commutes with the action of $\sl(2,\C)$, $V_{\kappa\tau}$ is also a
representation.

Recall that $Z_{\kappa\tau}(a,a,2b)$ is constructed by deleting vertices
fixed by $\kappa\tau$ and edges incident to those vertices, and then
taking one of the two isomorphic components of the complement.  By
inspection of the basis for $V_{\kappa\tau}$, the restriction
$\alpha_{\kappa\tau}(X)|_{-1}$ of
$\alpha(X)|_{-1}$ to $V_{\kappa\tau}$ is a weighted, bipartite adjacency
matrix for $Z_{\kappa\tau}(a,a,2b)$.  As in the no-symmetry case, opposite
edges of each hexagonal face have the same weight, and therefore it
is Kasteleyn-flat.  Indeed, the weighting of $Z_{\kappa\tau}(a,a,2b)$ is
induced from that of $Z(a,a,2b)$.  It follows that
\bqn N_{\kappa\tau}(a,a,2b) =
\det(\alpha_{\kappa\tau}(X)|_{-1})/m_{\kappa\tau},\label{enkappatau}\eqn
where $m_{\kappa\tau}$ is the value of a single term in the determinant
(\ie, the weight of any matching in $Z_{\kappa\tau}(a,a,2b)$).

To compute $m_{\kappa\tau}$, observe that a $\kappa\tau$-invariant
matching of $Z(a,a,2b)$ has two edges of the same weight for each edge
in the corresponding
matching of $Z_{\kappa\tau}(a,a,2b)$, plus all of the deleted
edges.  The weights of the deleted edges are $1,3,5,\ldots,2a-1$.
Therefore
\bqn m_{\kappa\tau} = \sqrt{|m|/(2a-1)!!},\label{etkappatau}\eqn
where $n!! = (n-1)(n-3)\ldots$ is the odd factorial.  The
determinant of $\alpha(X)|_{-1}$ is best computed using the $D$
map of Section~\ref{ssl2}.  The representation $V$ has two
basis vectors $e_i \tensor e_j \tensor e_k$ for each basis vector of
$V_{\kappa\tau}$, plus  basis vectors of the form
$e_i \tensor e_i \tensor e_k$.
Therefore
$$\ch(V_{\kappa\tau}) = \frac12\left(\ch(V)-\sum_{i=1}^{a+2b}
\sum_{k=1}^{2a} t^{4i+2k-4a-4b-3}\right).$$
Let $P(t) = \ch(V) - 2\ch(V_{\kappa\tau})$ be the double summation; then
\begin{align*}
D(P(t)) &= \prod_{i=1}^{a+2b}\prod_{k=1}^{2a}
    \frac{\max(2i+k-2a-2b-1,1)^2}{\max(2i+k-1-2a-2b-2,1)^2} \\
&= \prod_{i=1}^{a+2b} \frac{\max(2i-2b-1,1)^2}{\max(2i-2a-2b-2,1)^2} \\
&= \frac{(2a+2b-1)!!^2}{4^{b-1}(b-1)!^2}.
\end{align*}
Therefore
$$D(\ch(V_{\kappa\tau})) =\sqrt{\frac{4^{b-1}(b-1)!^2D(\ch(V))}{(2a+2b-1)!!^2}}.$$
Combining with equations~(\ref{enkappatau}) and (\ref{etkappatau}),
$$N_{\kappa\tau}(a,a,2b) =
\sqrt{\frac{N(a,a,2b) 2^{b-1}(b-1)!(2a-1)!!}{(2a+2b-1)!!}}$$
This is equivalent to standard formulas for $N_{\kappa\tau}(a,a,2b)$
\cite{Stanley:symmetries}.

Similarly, the symmetry $\rho$ on $Z(a,a,a)$ is induced by the linear
transformation
$$L_\rho \in \End(V_{2a-1} \tensor V_{2a-1} \tensor V_{2a-1})$$
given by 
$$L_{\rho}(e_i \tensor e_j \tensor e_k) = e_j \tensor e_k \tensor e_i.$$
Define $V_\rho$ to be the eigenspace of $L_{\rho}$ with eigenvalue $1$,
and choose as a basis for $V_\rho$ vectors of the form $e_i \tensor e_i 
\tensor e_i$ and vectors of the form
$$e_i \tensor e_j \tensor e_k + e_j \tensor e_k \tensor e_i +
e_k \tensor e_i \tensor e_j$$
when $i$, $j$, and $k$ are not all equal. Since $Z_\rho(a,a,a)$ is the
quotient graph $Z(a,a,a)/\rho$, the restriction $\alpha_\rho(X)|_{-1}$
of $\alpha(X)|_{-1}$ to
$V_\rho$ must be a weighted, bipartite adjacency matrix for $Z_\rho(a,a,a)$. 
Again, the weighting of $Z_\rho(a,a,a)$ induced from $Z(a,a,a)$ is compatible
with $\alpha_\rho(X)|_{-1}$ and is Kasteleyn-flat.  If  $m_\rho$ is
the weight of any matching in $Z_{\rho}(a,a,a)$, then
$$m_\rho = \sqrt[3]{m}.$$
Meanwhile,
$$D(\ch(V_\rho)) = \sqrt[3]{D(\ch(V))
\prod_{i=1}^{2a} \frac{\max(6i-3a-1,1)^4}{\max(6i-3a-2,1)^4}},$$
because
$$\ch(V_\rho)(t) = \frac13\left(\ch(V)(t) + 2\ch(V_{2a-1})(t^3)\right)$$
by counting basis vectors.  The final result is
$$N_\rho(a,a,a) = \sqrt[3]{N(a,a,a)\prod_{i=1}^a \frac{(3i-1)^2}{(3i-2)^2}}.$$
This is again easily equivalent to the standard formulas
\cite{Stanley:symmetries}.

In $V = V_{4a-1}^{\tensor 3}$, we define
$$V_{\kappa\tau,\rho} = V_{\kappa\tau} \cap V_\rho = \Lambda^3 V_{2a-1},$$
together with the preferred basis of vectors $e_i \wedge e_j \wedge
e_k$ with $i < j < k$.  (Recall that the wedge product is the
anti-symmetrized tensor product.) In this basis,
$\alpha_{\kappa\tau,\rho}(X)|_{-1}$ is a Kasteleyn-flat matrix for
$Z_{\rho,\kappa\tau}(2a,2a,2a)$.  A product formula for
$N_{\rho,\kappa\tau}(2a,2a,2a)$ is the inevitable consequence, and in
this final case, we omit the details.


\providecommand{\bysame}{\leavevmode\hbox to3em{\hrulefill}\thinspace}

\end{document}